\def\vers{Apr.~8, 2009, v.4}
\magnification=1200
\hsize=6.5truein
\vsize=8.9truein
\font\bigfont=cmr10 at 14pt
\font\mfont=cmr9
\font\sfont=cmr8
\font\mbfont=cmbx9
\font\sifont=cmti8
\def\scirc{\,\raise.2ex\hbox{${\scriptstyle\circ}$}\,}
\def\mopl{\hbox{$\bigoplus$}}
\def\msum{\hbox{$\sum$}}
\def\a{\alpha}
\def\C{{\bf C}}
\def\D{\Delta}

\def\DD{\widetilde{\Delta}}
\def\e{\varepsilon}
\def\f{\bar{f}}
\def\g{\bar{g}}
\def\G{\Gamma}
\def\GG{\widehat{\Gamma}}
\def\H{{\bf H}}
\def\h{\eta}
\def\i{\iota}
\def\L{{\cal L}}
\def\LL{\widehat{\cal L}}
\def\m{\mu}
\def\mm{\widetilde{\mu}}
\def\n{\nu}
\def\nn{\widetilde{\nu}}
\def\O{{\cal O}}
\def\Q{{\bf Q}}
\def\R{{\bf R}}
\def\S{\bar{S}}
\def\SS{\widetilde{S}}
\def\s{\sigma}
\def\V{{\cal V}}
\def\VV{\widehat{\cal V}}
\def\Z{{\bf Z}}
\def\z{\bar{z}}

\def\Ext{\hbox{\rm Ext}}
\def\Hom{\hbox{\rm Hom}}
\def\Ker{\hbox{\rm Ker}}
\def\Gr{\hbox{\rm Gr}}
\def\GGK{{\rm GGK}}
\def\Inv{{\rm Inv}}

\def\Re{{\rm Re}}
\def\Im{{\rm Im}}
\def\1{{\hskip1pt}}
\def\simto{\buildrel\sim\over\longrightarrow}
\def\onto{\to\!\!\!\!\!\to}
\hbox{}
\vskip 1cm

\centerline{\bigfont Hausdorff property of the N\'eron models}

\smallskip
\centerline{\bigfont of Green, Griffiths and Kerr}

\bigskip
\centerline{Morihiko Saito}

\bigskip\medskip
{\narrower\noindent
{\mbfont Abstract.} {\mfont
We prove the Hausdorff property of the N\'eron modle of the family
of intermediate Jacobians which is recently defined by Green,
Griffiths and Kerr assuming that the divisor at infinity is smooth.
Using their result, this implies in this case the analyticity of
the closure of the zero locus of an admissible normal function.
The last assertion is also obtained by Brosnan and Pearlstein
generalizing their method in the curve case.}
\par}

\bigskip\bigskip
\centerline{\bf Introduction}
\footnote{}{{\sifont Date\1}{\sfont:\ \vers}}

\bigskip\noindent
Let $\H$ be a polarizable variation of $\Z$-Hodge structure of
weight $w<0$ on a complex manifold $S$.
Let $J_S(\H)$ denote the family of intermediate Jacobians over $S$.
A {\it normal function} is a holomorphic section of $J_S(\H)$
satisfying the Griffiths transversality, and corresponds by a
well-known theorem of Carlson [3] to an extension class of the
trivial variation $\Z_S$ by $\H$ in the abelian category of mixed
$\Z$-Hodge structures [6].
A normal function is called {\it admissible} [10] if the
corresponding extension is an admissible variation of mixed Hodge
structure in the sense of [9], [13].
This condition corresponds to certain conditions on the local
monodromy and the logarithmic growth at infinity of the normal
function, see [8].
By M.~Green and P.~Griffiths (see [1]) we have the following

\medskip\noindent
{\bf Conjecture~1.} The zero locus of an admissible normal function
is algebraic if $S$ is algebraic (where $w=-1$).

\medskip
Let $\S$ be a smooth partial compactification of $S$ such that
$D:=\S\setminus S$ is a divisor with normal crossings.
Then Conjecture~1 is easily reduced to

\medskip\noindent
{\bf Conjecture~2.} The closure of the zero locus of an admissible
normal function in $\S$ is analytic (where $w=-1$).

\medskip
In the curve case these conjectures were solved by
P.~Brosnan and G.~Pearlstein [1] where the assertion was equivalent
to the discreteness of the closure.
In this paper we prove

\medskip\noindent
{\bf Theorem~1.} {\it
Assume $D$ is smooth. If $w=-1$, then Conjecture~$2$ is true, and
moreover, the closure of the zero locus is contained in $S$ if the
local cohomological invariant of the admissible normal function is
nonzero and the local monodromy of $\H$ is unipotent.
The same assertions hold for $w<-1$ if we assume further that the
local cohomological invariant of the admissible normal function is
torsion.}

\medskip
We are informed that the first assertion is also obtained by
P.~Brosnan and G.~Pearlstein generalizing their method in [1],
see [2].
From Theorem~1 we can deduce

\medskip\noindent
{\bf Corollary~1.} {\it Assume $S$ admits a smooth compactification
$\S$ such that $D=\S\setminus S$ is smooth.
Then Conjecture~$1$ holds where $w=-1$.
For $w<-1$ it holds if furthermore the local cohomological invariant
of the admissible normal function is torsion.}

\medskip
The proof of Theorem~1 is reduced to the case where the local
monodromy is unipotent.
Let $J^Z_{\S}(\H)$ denote the Zucker extension of $J_S(\H)$, see [14].
Let $J^{\GGK}_{\S}(\H)^0$ denote the subspace of $J^Z_{\S}(\H)|_D$
defined by the quotient of the monodromy invariant subspace.
This is the identify component of the N\'eron model
$J^{\GGK}_{\S}(\H)$ which is recently defined by M.~Green,
P.~Griffiths and M.~Kerr, see [8].
Using their result, Theorem~1 is then reduced to the
following

\medskip\noindent
{\bf Theorem~2.} {\it
The Zucker extension $J^Z_{\S}(\H)$ has a structure of a complex Lie
group over $\S$, and it is a Hausdorff topological space on a
neighborhood of $J^{\GGK}_{\S}(\H)^0$.}

\medskip
This is essentially equivalent to Corollary~(2.7) which is
deduced from Theorem~(2.1).
The latter is proved by generalizing some arguments in [10]
used to define the structure of a complex Lie group on the Zucker
extension $J^Z_{\S}(\H)$ in the curve case.
The proof consists of elementary calculations for norm estimates.
Note that the Zucker extension outside the monodromy invariant
subspace is not necessarily Hausdorff, see e.g. [10], Remark~3.5~(iv).
As a corollary we get

\medskip\noindent
{\bf Theorem~3.} {\it
The N\'eron model of Green, Griffiths and Kerr $J^{\GGK}_{\S}(\H)$
is Hausdorff where $w=-1$.}

\medskip

In Section~1 we review some basic facts from the theory of Zucker
extensions, admissible normal functions and limit mixed Hodge
structures, and prepare some notation.
In Section~2 we first show Theorem~(2.1) on the norm estimate, and
then prove Theorems~1--3.

\bigskip\bigskip
\centerline{\bf 1.\ Preliminaries}

\bigskip\noindent
In this section we review some basic facts from the theories of
Zucker extensions, admissible normal functions, and limit mixed
Hodge structures, and prepare some notation.

\medskip\noindent
{\bf 1.1.~Family of intermediate Jacobians.}
Let $S$ be a complex manifold of dimension $n$, and
let $\H=((\L,F),L_{\Z})$ be a polarizable variation of $\Z$-Hodge
structure of weight $w<0$ on $S$. Here $(\L,F)$ is the underlying
filtered locally free sheaf and $L_{\Z}$ is the underlying $\Z$-local
system.
We assume that $L_{\Z}$ is torsion-free in this paper.
Let $\V$ be the vector bundle on $S$ corresponding to the locally
free sheaf $\L/F^0\L$, and let $\G\subset\V$ denote the subgroup
over $S$ corresponding to the subsheaf $L_{\Z}\subset\L/F^0\L$.
where the last injectivity follows from the condition that the
weight of $\H$ is negative.
The family of intermediate Jacobians is defined by
$$J_S(\H)=\V/\G,$$
where the quotient is set-theoretically taken fiberwise.
This has a structure of a complex manifold over $S$.

\medskip\noindent
{\bf 1.2.~Zucker extensions.}
Let $\S$ be a partial compactification of $S$ such that $\S$ and
$D:=\S\setminus S$ are smooth.
Let $j:S\to\S$ denote the inclusion.
Let $\LL$ be the Deligne extension of $\L$ such that the eigenvalues
of the residues of the logarithmic connection are contained in $[0,1)$,
see [5].
By Schmid [11], the Hodge filtration $F$ on $\L$ is uniquely extended
to a filtration $F$ on $\LL$ such that $\Gr_F^p\LL$ are locally free
(i.e. $F^p\LL=\LL\cap j_*F^p\L$).
Let $\pi:\VV\to\S$ denote the vector bundle corresponding to the
locally free sheaf $\LL/F^0\LL$.
Let $\GG$ be the subgroup of $\VV$ over $\S$ corresponding to
$j_*L_{\Z}\subset\LL/F^0\LL$.
We define the Zucker extension by
$$J_{\S}^Z(\H)=\VV/\GG.$$
We will show in Section 2 that it has a structure of a complex Lie
group over $\S$, see [10] for the curve case.

Let $\VV_D$ be the restriction of $\VV$ over $D$, and $\VV_D^{\Inv}$
be the union of the subspaces defined by $\Ker\,N\subset \LL(s):=
\LL/m_s\LL$ for $s\in D$, where $m_s$ is the maximal ideal at $s$.
The subspace of $J_{\S}^Z(\H)$ defined by the quotient of
$\VV_D^{\Inv}$ will be denoted by $J^{\GGK}_{\S}(\H)^0$.

\medskip\noindent
{\bf 1.3.~Admissible normal functions.}
A normal function is a holomorphic section $\n$ of $J_S(\H)$ satisfying
the Griffiths transversality.
By Carlson [3] it corresponds to an extension class of $\Z_S$ by $\H$
giving a short exact sequence
$$0\to\H\to\H'\to\Z_S\to 0.\leqno(1.3.1)$$
A normal function $\n$ is called admissible with respect to a partial
compactification $\S$ of $S$ if $\H'$ is an admissible variation of 
mixed Hodge structure ([9], [13]) with respect to $\S$, see [10].
Here $\S$ is as in (1.2).
In the unipotent monodromy case $\H'$ should satisfy the following
two conditions:

\smallskip
(i) $\Gr_F^p\Gr_k^W\LL'$ are locally free.

\smallskip
(ii) The relative monodromy filtration exists.

\smallskip\noindent
By [8] these conditions in case $w=-1$ are equivalent respectively to

\smallskip
(i)$'$ A lifting $\nn$ of $\n$ in $\V$ has logarithmic growth.

\smallskip
(ii)$'$ The local cohomological invariant $\gamma(\n)$ is torsion,
see also [10], Remark~1.6 (iv).

\smallskip\noindent
Here the local cohomological invariant $\gamma(\n)\in H^1(S,L_{\Z})
=\Ext^1(\Z_S,L_{\Z})$ is defined by passing to the underlying short
exact sequence of $\Z$-local systems of (1.3.1) after shrinking $\S$
sufficiently, or by considering the `monodromy' of a lifting $\nn$ of
$\n$ in $\V$, see [8].

It is known that the normal function is extended to a section of
the Zucker extension if and only if the local cohomological invariant
$\gamma(\n)$ vanishes, see e.g. [10], Prop.~2.3 and 2.4
(and [7] for the geometric case).
For example, if the local cohomological invariant vanishes, then
the normal function is given by the difference between two splittings
$\s_{\Z},\s_F$ of the underlying short exact sequence of locally free
sheaves of (1.3.1)
$$0\to\LL\to\LL'\to\O_{\S}\to 0,$$
where $\s_{\Z}$ is defined over $\Z$ and $\s_F$ is compatible with
$F$. Moreover we have

\medskip\noindent
{\bf Theorem} (Green, Griffiths, Kerr [8]).
The above normal function passes through the subspace
$J^{\GGK}_{\S}(\H)^0$.

\medskip
This is proved by using the Griffiths transversality,
see [8] for details.

\medskip\noindent
{\bf 1.4.~Limit mixed Hodge structure.}
From now on, assume $S=\D^*\times\D^{n-1}$ and $\S=\D^n$ where $\D$
is a polydisk.
Let $t=(t_1,t_2,\dots,t_n)$ be the coordinate system of $\D^n$.
Set $S'=\D^{n-1}$ and $t'=(t_2,\dots,t_n)$.

Let $\H=((\L,F),L_{\Z})$ be a polarizable variation of $\Z$-Hodge
structure of weight $w<0$, and let $(\LL,F)$ be as in (1.2).
Assume the local monodromy $T$ around $t_1=0$ is unipotent.
Set $N=\log T$.
By Schmid [11] we have the limit mixed Hodge structure of $\H$ at $0$
$$H=(H_{\Z},(H_{\Q},W),(H_{\C};W,F))\,\,\,(\hbox{or $(H_{\Z},W,F)$
to simplify the notation}).$$
By [5] (and [12] in the geometric case) $H_{\C}$ can be defined as
$$H_{\C}=\LL(0)\,(:=\LL/m_0\LL),$$
where $m_0$ is the maximal ideal at $0$,
and $F$ is the quotient filtration of $F$ on $\LL$.

For $A=\Z,\Q,\R$, $H_A$ is defined by
$$H_A=\G(\DD^*\times S',(\rho\times id)^*L_A),
\leqno(1.4.1)$$ where
$$\rho:\DD^*:=\{z\in\C\,|\,\Re\,z>r_0\}\to\D^*$$ is a universal
covering defined by $z\mapsto\exp(2\pi iz)$ (with $r_0>0$) and
$L_A=L_{\Z}\otimes_{\Z}A$ for $A=\Q,\R$.
We have the isomorphism (see [5])
$$H_{\C}=H_A\otimes_A\C,$$ induced by the injection
$H_A\to H_{\C}$ associating
$$\tilde{u}=\exp(-zN)u\in\LL\quad\hbox{to}\,\,\,u\in H_A,
\leqno(1.4.2)$$
where $z={1\over 2\pi i}\log t_1$.
Note that (1.4.2) gives a trivialization of the locally free sheaf
$\LL$ (which depends on the coordinate $t_1$).
In particular, $\{\tilde{u}_i\}$ is a basis of $\LL$ if $\{u_i\}$
is a basis of $H_A$.

The weight filtration $W$ is given by the monodromy filtration
associated to $N=\log T$ up to the shift by $w$.
This is characterized by the two conditions:
$N(W_jH_A)\subset W_{j-2}H_A$ for any $j$, and
$$N^i:\Gr^W_{w+i}H_A\simto\Gr^W_{w-i}H_A\quad\hbox{for}
\,\,\,i>0.
\leqno(1.4.3)$$
If $A=\C$, this is compatible with $F$ (up to a shift).
We may assume for the proof of Theorem~(2.1) that (1.4.3) holds also
for $A=\Z$ replacing $L_{\Z}$ if necessary (using the fact that
the action of $N$ on $\Gr^W$ coincides with that of $T-id$).
Note that
$$N:(H_{\Z},W,F)\to (H_{\Z},W,F)(-1)$$ is a morphism of mixed
$\Z$-Hodge structures (see [6] for the Tate twist), and hence
$K:=\Ker\,N$ is a mixed $\Z$-Hodge structure of weights $\le w$.
Since $w<0$, we have the injectivity of
$$K_{\Z}\to K_{\C}/F^0K_{\C}.$$

\medskip\noindent
{\bf 1.5.~Two splittings of the weight filtration.}
For $A=\Z,\Q,\R,\C$, set
$$G_A:=\mopl_i\Gr_i^WH_A\quad\hbox{with}\quad
W_kG_A:=\mopl_{i\le k}\Gr_i^WH_A.
\leqno(1.5.1)$$
If $A=\C$, it has the induced filtration $F$.
We choose and fix isomorphisms compatible with the action of $N$
$$ \a_{\R}:(G_{\R},W)\simto (H_{\R},W),\quad
\a_{\C}:(G_{\C};W,F)\simto (H_{\C};W,F),
\leqno(1.5.2)$$
inducing the identity on the graded pieces of $W$ and such that
$\a_{\R}$ is defined over $\Z$.
The existence of $\a_{\C}$ is equivalent to that of a basis
$(v_{i,p,k})$ of $\Gr_F^pP\Gr^W_{w+i}H_{\C}$ such that
$v_{i,p,k}$ is represented by an element of $F^pW_{w+i}H_{\C}$
annihilated by $N^{i+1}$ where $P\Gr^W_{w+i}$ denotes the primitive
part defined by $\Ker\,N^{i+1}\subset\Gr^W_{w+i}$.
Here the problem is the annihilation by $N^{i+1}$,
and the existence is reduced to the surjections
$$N^{i+1}:F^pW_{w+i-1}H_{\C}\onto F^{p-i-1}W_{w-i-3}H_{\C}
\quad\hbox{for}\,\,\,i>0.$$ 
This is further reduced to the surjections
$$N^{i+1}:F^i\Gr^W_{w+j}H_{\C}\onto F^{p-i-1}\Gr^W_{w+j-2i-2}H_{\C}
\quad\hbox{for}\,\,\,j\le i-1,$$ and follows from (1.4.3).
The argument is similar for $\a_{\Q}$, and it induces $\a_{\R}$.
We may assume for the proof of Theorem~(2.1) that $\a_{\Q}$ induces
the isomorphism $\a_{\Z}$ replacing $L_{\Z}$ if necessary.
However, $\a_{\R}$ is not compatible with $\a_{\C}$ unless $W$ splits
over $A$ in a compatible way with the Hodge filtration $F$.
The incompatibility is expressed by the injective morphism
$$\i:=\a_{\C}\scirc\a_A^{-1}:(G_A,W)\to(G_{\C},W)\quad(A=\Z,\Q,\R).
\leqno(1.5.3)$$

\medskip\noindent
{\bf 1.6.~Primitive decomposition.}
For $A=\Q,\R,\C$, we have the primitive decomposition
$$\Gr_i^WH_A=\mopl_{k\ge 0}\,N^kP\Gr_{i+2k}^WH_A,
\leqno(1.6.1)$$
where $P\Gr_{w+i}^WH_A$ is the primitive part defined by
$\Ker\,N^{i+1}\subset\Gr_{w+i}^WH_A$, and $w$ is the weight of $\H$.
We may assume for the proof of Theorem~(2.1) that (1.6.1) holds
also for $A=\Z$ replacing $L_{\Z}$ if necessary.
With the notation of (1.5.1), set
$$G_{A,k}^{(j)}:=N^{j-k}P\Gr_{w+j}^WH_A,\quad
G_A^{(j)}:=\mopl_{0\le k\le j}\,G_{A,k}^{(j)},\quad
G_{A,k}:=\mopl_{k\le j\le m}\,G_{A,k}^{(j)}.$$
Then we have the bigrading
$$G_A=\mopl_{0\le j\le m}\,G_A^{(j)}=\mopl_{0\le k\le m}\,G_{A,k}=
\mopl_{0\le k\le j\le m}\,G_{A,k}^{(j)}.$$
This is compatible with $F$ if $A=\C $.
Note that the $G_{A,k}$ give a splitting of the kernel filtration
$K_k$ on $G_A$ defined by $\Ker\,N^{k+1}$.

Since $\i$ in (1.5.3) is compatible with the action of $N$,
there are morphisms
$$\i_{i,j}: G_A^{(j)}\to G_{\C}^{(i)},$$
such that the restriction of $\i$ to $G_A^{(j)}$ is $\sum_i\i_{i,j}$.
Furthermore we have
$$\i_{i,j}=\msum_k\,\i_{i,j,k}\quad \hbox{with}\quad
\i_{i,j,k}G_{A,k'}^{(j)}\subset G_{A,k'-k}^{(i)}.$$
Here $\i_{i,j,k}$ are compatible with the action of $N$, and
$$\i_{i,j,k}=0\quad\hbox{if}\,\,\,k<0.
\leqno(1.6.2)$$
This is equivalent to that $\i$ preserves the kernel filtration
$K_k$ defined by $\Ker\,N^{k+1}$.

\medskip\noindent
{\bf 1.7.~Hodge decomposition.}
We have the Hodge decomposition
$$G_{\C}=\mopl_{i,p}\,(\Gr_i^WH_{\C})^{p,i-p}.
\leqno(1.7.1)$$
This is compatible with the action of $N$ on $G_{\C}$. Set
$$G_{\C}^{<0}:=\mopl_{i\in\Z,p<0}\,(\Gr_i^WH_{\C})^{p,i-p}\subset
G_{\C}.$$
Since the Hodge decomposition is compatible with the primitive
decomposition (1.6.1), we have
$$G_{\C}^{<0}=\mopl_{0\le j\le m}\,G_{\C}^{<0,(j)}=
\mopl_{0\le k\le j\le m}\,G_{\C,k}^{<0,(j)}$$
with $G_{\C,k}^{<0,(j)}=G_{\C}^{<0}\cap G_{\C,k}^{(j)},$ etc.

\medskip\noindent
{\bf 1.8.~Trivialization of the vector bundle.}
Identifying $H_{\C}$ with $G_{\C}$ by $\a_{\C}$ in (1.5.2) and
using the trivialization of $\LL$ induced by (1.4.2), the subspace
$G_{\C}^{<0}$ of $G_{\C}$ in (1.7) defines a trivialized subsheaf
$\LL^{<0}$ of $\LL$ such that
$$\LL=\LL^{<0}\oplus F^{0}\LL,
\leqno(1.8.1)$$
shrinking $\D$ if necessary.
This is equivalent to the trivialization
$$\VV=G_{\C}^{<0}\times\S .
\leqno(1.8.2)$$
By (1.5.2), an element $v$ of $G_{\R}$ is identified with an element
of $H_{\R}$, and determines a multivalued section of
$\V=\VV|_S$.
Using the trivialization (1.8.2), it is expressed as
$$\phi(v)\in G_{\C}^{<0}\otimes_{\C}\G(\SS,\O_{\SS}),
\leqno(1.8.3)$$
where $\SS=\DD^*\times S'$ in the notation of (1.4.1).
For $(z,t')\in\SS$, the value of $\phi(v)$ at $(z,t')$ will be
denoted by $\phi(v;z,t')$.

\medskip\noindent
{\bf 1.9.~Norms.}
We choose and fix a norm on each complex vector space
$G_{\C,k}^{<0,(j)}$, i.e. there is a continuous map
$v\mapsto|v|\in\R_{\ge 0}$ satisfying the usual relations:
$|v+v'|\le|v|+|v'|$, $|av|=|a||v|$ for $a\in\C$, and $\{|v|=1\}$
is compact.
We may assume that it is compatible with the Hodge decomposition
(1.7.1) (i.e. the norm is the sum of the norm of the direct factors).
This induces a norm on $G_{\C}^{<0}$ compatible with the
primitive and Hodge decompositions.

For $v=\sum_{0\le k\le j\le m}v_k^{(j)}\in G_{\C}^{<0}$ with
$v_k^{(j)}\in G_{\C,k}^{<0,(j)}$, set
$$|v|_k:=|v_k|=\msum_{0\le j\le m}\,|v_k^{(j)}|\quad\hbox{with}\,\,\,
v_k=\msum_{k\le j\le m}\,v_k^{(j)}.\leqno(1.9.1)$$

We also choose and fix a norm on each real vector space
$P\Gr_{w+i}^WH_{\R}$.
This induces a norm on $G_A$ using the primitive decomposition
(1.6.1). Let
$$u=\msum_{0\le k\le i\le m}\,u^{(i)}_k
=\msum_{0\le k\le m}\,u_k\in G_{\R},$$
where $u^{(i)}_k=N^{i-k}u_{i,k}$ with $u_{i,k}\in P\Gr_{w+i}^WH_{\R}$
and $u_k=\sum_{k\le i\le m}u^{(i)}_k$.
Let $z\in\DD^*$, and set
$$\eqalign{E(G_{\Z})&:=\min\{|u|:u\in G_{\Z}\setminus\{0\}\},
\quad y:=\Im\,z,\cr
A(u,z)&:=\msum_{0\le k\le m}\,|u_{k}|y^k=\msum_{0\le k\le i\le m}\,
|u_{i,k}|y^k,\cr
B(u,z)&:=\msum_{1\le k\le m}\,|u_{k}|y^{k-1}=
\msum_{1\le k\le i\le m}\,|u_{i,k}|y^{k-1}.}
\leqno(1.9.2)
$$
For $\h>0$ and $r>1$, set
$$I_{r,\h}=\{(z,t')\in\C\times S'\,|\,0\le\Re\,z<1,\Im\,z>r,|t_j|
<\h\,(j>1)\}\subset\SS.\leqno(1.9.3)$$
In the case $S'=pt$, we will denote $I_{r,\h}$ by $I_r$.

\medskip\noindent
{\bf 1.10.~Examples.}
(i) Assume $H_{\Z}=\Z e_0+\Z e_1$ with $Ne_0=e_1$,
$Ne_1=0$, and $F^0=\O_{\D}\tilde{e}_0\,$ (i.e. an $\R$-split
nilpotent orbit). Then
$\VV=\C\times\D$, and $e_0,e_1$ are identified respectively with the (multivalued) functions $1$ and $z:={1\over 2\pi i}\log t$ on $\D^*$.
(This is the simplest example with multivalued 
\medskip\noindent
{\bf Examples~2.}
$H_{\Z}=\Z e_0+\Z e_1+\Z f_0+\Z f_1\,\,$ s.t.

\smallskip
$Ne_0=e_1$, $Ne_1=0$, $Nf_0=f_1$, $Nf_1=0$,

\smallskip
$F^0=\O_{\D}\tilde{v}_0+\O_{\D}\tilde{v}_1\,$ with
$\,v_k=e_k-if_k\,\,(k=0,1)$.

\smallskip
Then $\,\VV=\C^2\times\D\,\,$ s.t. $\,\G\,\Leftrightarrow\,
\{(az+b,a)\,|\,a,b\in\Z[i]\}$.

\smallskip
So $\VV/\GG$ is $\underline{\hbox{not}}$ Hausdorff.

Indeed, $(\gamma,a)\in\overline{\G}\subset\VV$
for any $\gamma\in\C$, $a\in\Z\setminus\{0\}$,
considering $(az-ni,a)$
for $z=(\gamma+ni)/a\,\,(n\gg 1)$.

\bigskip\bigskip
\centerline{\bf 2.\ Norm estimates}

\bigskip\noindent
In this Section we prove Theorems~1 and 2.
We first show the following which is a key to the proof of Theorem~2.

\medskip\noindent
{\bf 2.1.~Theorem.} {\it With the notation of {\rm (1.5--9)}, let
$v\in G_{\C}^{<0}\cap\Ker\,N$. Then there exist $\e,\h>0$ and $r>1$
such that for any $u=\sum_{j,k}{u}_k^{(j)}\in G_{\Z}\setminus\Ker\,N$
and for $(z,t')\in I_{r,\h}$, the following holds\1$:$

\smallskip\noindent
{\rm (i)} We have for some $k\ge 0$
$$|\phi(\i(u);z,t')-v|_k\ge\e A(u,z)/y^k.\leqno(2.1.1)$$

\smallskip\noindent
{\rm (ii)} If $(2.1.1)$ holds for some $k>0$,
then $u_{k'}^{(j)}\ne 0$ for some $k'\ge k$.}

\medskip\noindent
{\bf 2.2.~Remarks.} (i) Condition (2.1.1) is essentially
independent of the choice of norms on $G_{\C,k}^{<0,(j)}$ and
$P\Gr_{w+i}^WH_{\R}$, replacing $\e$ if necessary.

\medskip\noindent
(ii) If $u\in\Ker\,N$ and (2.1.1) does not hold for every $k$, then
$$\i(u)-v\in F^0G_{\C},$$
assuming that $\e$ is sufficiently small compared to the distance
between the images of $v$ and $\i(G_{\Z}\cap\Ker\,N)\setminus\{v\}$
in $G_{\C}/F^0G_{\C}$.

\medskip\noindent
{\bf 2.3~Proof of Theorem~(2.1)(i).}
We start with the following

\medskip\noindent
{\bf (a) Reduction to the case of a nilpotent orbit with $S'=pt$.}

\smallskip\noindent
Using the trivialization of the vector bundle $\LL$ given by (1.4.2),
there is
$$M(t)\in\Hom_{\C}(F^{0}G_{\C}, G_{\C}^{<0})\otimes_{\C}
\G(\S,\O_{\S}),$$
such that the fiber of $F^{0}\LL$ at $t\in\S$ is expressed by the
graph of $M(t)$ using the following direct sum decomposition in the
notation of (1.7):
$$G_{\C}=G_{\C}^{<0}\oplus F^0G_{\C}.$$
Note that $M(0)=0$.
Here $G_{\C}$ is identified with $H_{\C}$ by $\a_{\C}$ in (1.5.2).

By the inverse of (1.4.2), $\tilde{v}\in G_{\C}$ corresponds to a
multivalued section $ v$ of $\L $, which is expressed as
$$\msum_{k\ge 0}\,N^k\tilde{v}\otimes z^k/k!\in G_{\C}\otimes_{\C}
\G(\SS,\O_{\SS}).\leqno(2.3.1)$$
Using the above decomposition, its image in $\VV$ (see (1.8.3)) is
given by
$$\eqalign{\phi(v)&=\phi'(v)-M(t)\phi''(v)\in G_{\C}^{<0}\otimes_{\C}
\G(\SS,\O_{\SS})\quad\hbox{with}\cr
\phi'(v)&=\msum_{k\ge 0}\,pr_1(N^k\tilde{v})\otimes z^k/k!,\quad
\phi''(v)=\msum_{k\ge 0}\,pr_2(N^k\tilde{v})\otimes z^k/k!,}$$
where $pr_1:G_{\C}\to G_{\C}^{<0}$ and $pr_2:G_{\C}\to F^{0}G_{\C}$
are the projections associated to the above decomposition.
By (1.6.2) we have
$$|\phi''(\i(u);z,t')|_k\le CA(u,z)/y^k\le CA(u,z)\quad\hbox{for any}
\,\,\,0\le k\le m,$$
where $C>0$ is independent of $u,z,t'$
(and $1\le y\le|z|\le 2y$ because $r\ge 1$).
Since $M(0)=0$, this implies
$$|M(t)\phi''(\i(u);z,t')|_k\le C'(\msum_{j\ge 0}\,|t_j|)A(u,z)\quad
\hbox{for any}\,\,\,0\le k\le m,\leqno(2.3.2)$$
where $C'>0$ depends only on $C,M(t)$, and is independent of $u,z,t'$.
Then
$$|M(t)\phi''(\i(u);z,t')|_k\le\e A(u,z)/y^k\quad\hbox{for any}\,\,\,
0\le k\le m,$$
if $1/r$ and $\h$ are sufficiently small (depending only on $C',m$
and $\e$).
So we may replace $\phi$ and $\e$ in (2.1.1) by
$\phi'$ and $2\e$ respectively.
Thus the assertion is reduced to the case $\H$ is a nilpotent orbit,
where $M(t)=0$ and $\phi''=0$.
(Note that a variation of Hodge structure is a nilpotent orbit
if and only if there is a basis $\{v_i\}$ of $H_{\Z}\otimes\C$
such that the $\tilde{v}_i$ in (1.4.2) give generators of
$\bigoplus_p\Gr_F^p\LL$.)
Then we may assume moreover
$$S'=pt,\quad S=\D^*.$$

From now on, $\phi(\i(u);z,t')$ will be denoted by $\phi(\i(u),z)$.

\medskip\noindent
{\bf (b) Reduction to the case:} $\i_{i,j,k}=0$ for $k\ne 0$.

\medskip\noindent
For any $u^{(j)}\in G_{\Z}^{(j)}$ and $z\in I_r$, we have
$$|\phi(\i_{i,j,k}(u^{(j)}),z)|_{k'}\le C''A(u^{(j)},z)/y^{k+k'},
\leqno(2.3.3)$$
where $C''>0$ is independent of $u^{(j)},z$.
So the assertion is reduced to the desired case by the same argument
as above.
In particular, we get
$$\i_{i,j}=0\quad\hbox{if}\,\,\,i<j,$$
and $\i_{i,j}$ are graded morphisms (up to the shift by $i-j)$.

\medskip\noindent
{\bf (c) Reduction to the case:} $G_{\R}=G_{\R}^{(i)}$.

\medskip\noindent
Let $v=\sum_iv^{(i)}$ be the decomposition such that $v^{(i)}\in
G_{\C}^{(i)}\cap\Ker\,N$.
Assume the assertion is true for each
$(G_{\R}^{(i)},G_{\C}^{(i)},\i_{i,i},v^{(i)})$.
Then there exist $\e',r'>0$ such that for any
$u^{(i)}\in G_{\Z}^{(i)}\setminus\Ker\,N$ and $z\in I_{r'}$, we have
$$|\phi(\i_{i,i}(u^{(i)}),z)-v^{(i)}|_{k(i)}\ge\e'A(u^{(i)},z)
/y^{k(i)}\quad\hbox{for some}\,\,k(i)\ge 0,\leqno(2.3.4)$$
where $k(i)$ may depend on $i$.
However, $\e',r'$ are independent of $i$
(taking the minimum of $\e$ and the maximum of $r$ for $i$).

Assume on the other hand the assertion is false for
$(G_{\R},G_{\C},\i,v)$. Then, for any $\e,r>0$, there exist
$u=\sum_iu^{(i)}\in G_{\Z}\setminus\Ker\,N$ and $z\in I_{r}$ such
that
$$|\phi(\i(u),z)-v|_k\le\e A(u,z)/y^k\quad\hbox{for any}\,\,\,
k\ge 0.$$
By (2.3.3) with $k=0$ and $k'$ replaced by $k$, we get
$$|\phi(\i_{i,j}(u^{(j)}),z)|_k\le C_{i,j}A(u^{(j)},z)/y^k
\quad\hbox{for any}\,\,\,i,j,k,$$
where $C_{i,j}>0$ is independent of $u^{(j)},z$.
We have furthermore for any $i,j,k$
$$|\phi(\i_{i,i}(u^{(i)}),z)-v^{(i)}|_k\le |\phi(\i(u),z)-v|_k+
\msum_{j<i}\,|\phi(\i_{i,j}(u^{(j)}),z)|_k.$$
In case (2.3.4) holds for $i$, we get thus
$$\e'A(u^{(i)},z)<\e A(u,z)+\msum_{j<i}\,C_{i,j}A(u^{(j)},z).
\leqno(2.3.5)$$
Here we do not necessarily assume $u^{(i)}\notin\Ker\,N$.
Note that (2.3.4) does not hold only in the case
$u^{(i)}\in\Ker\,N$ with $\i_{i,i}(u^{(i)})=v^{(i)}$ mod $F^0$
(assuming $\e'$ sufficiently small).
In this case we have $A(u^{(i)},z)=|u^{(i)}|$,
and this can be bounded essentially by $|v^{(i)}|$ using the
injectivity of
$$G_{\Z}\cap\Ker\,N\to G_{\C}/F^0.$$
Since $u\in G_{\Z}\setminus\Ker\,N$, we see that (2.3.5) holds
also in this case if we assume in the notation of (1.9.2) $r$ is
sufficiently large compared to
$$\e'|v^{(i)}|/\e E(G_{\Z}).$$
Thus (2.3.5) holds for any $i$.
This contradicts the following if $\e$ is sufficiently small,
where $a_i=A(u^{(i)},z)$, $C'_{i,j}=C_{i,j}/\e'$, and $\e''=\e/\e'$.

\medskip\noindent
{\bf Sublemma.} {\it
Let $C'_{i,j}$ be positive numbers for $0\le j<i\le m$.
Then there exists a positive number $\e''<1/2$ such that for any
nonnegative numbers $a_i\,(0\le i\le m)$,
the following condition implies $a_i=0:$}
$$a_i\le\e''(\msum_{j}\,a_{j})+\msum_{j<i}\,C'_{i,j}a_{j}\quad
\hbox{for}\,\,\,0\le i\le m.$$

\medskip\noindent
{\it Proof.}
The condition for $i=0$ implies
$$a_0\le 2\e''(\msum_{j>0}\,a_{j}),$$ and hence
$$a_i\le B\e''(\msum_{j>0}\,a_{j})+\msum_{0<j<i}\,C'_{i,j}a_{j}\quad
\hbox{for}\,\,\,1\le i\le m,$$
where $B>0$ depends only on the $C'_{i,0}$.
So the assertion follows by induction on $m$.

\medskip
Thus the assertion is reduced to the case $G_{\R}=G_{\R}^{(i)}$.

\medskip\noindent
{\bf (d) Reduction to the case:} $\dim P\Gr_{w+i}^WH_{\R}=1$ or $2$.

\medskip\noindent
We replace the assertion by the following which implies
Theorem~(2.1)(i) by the above arguments, since $u\in G_{\Z}$ belongs
to $\Ker\,N$ if we have in the notation of (1.9.2)
$$B(u,z)\le A(u,z)/r\le C_v/r<E(G_{\Z}).$$

\medskip\noindent
{\bf 2.4.~Proposition.} {\it
With the above notation, assume $\H$ is a nilpotent orbit on
$S=\D^*$, $G_{\R}=G_{\R}^{(i)}$ for some $i$, and $\i_{i,j,k}=0$
for $k\ne 0$. Then, for $v\in G_{\C}^{<0}\cap\Ker\,N$, there exist
$C_v,\e>0$ and $r>\max(C_v/E(G_{\Z}),1)$ such that for any
$u=\sum_{j,k}{u}_k^{(j)}\in G_{\R}$ and for $z\in I_{r}$, we have
$$A(u,z)\le C_v,$$ if the following condition is satisfied for any
$k\ge 0:$}
$$|\phi(\i(u),z)-v|_k\le\e A(u,z)/y^k.\leqno(2.4.1)$$

\medskip\noindent
{\it Proof.} Since the assertion depends only on the underlying
variation of $\R$-Hodge structure, we may assume that $\H$ is
defined only over $\R$ (forgetting the $\Q$-structure).
Furthermore, if $\H=\H'\oplus\H''$ and the assertion holds for
$\H',\H''$, then the assertion also holds for $\H$ (modifying $\e$
if necessary).
Indeed, for $u=u'+u''$ and $v=v'+v''$ with $u',v'$ coming from
$\H'$ and similarly for $u'',v''$, we have
$$|\phi(\i(u'),z)-v'|_k\le|\phi(\i(u),z)-v|_k\le 2\e A(u',z)/y^k,$$
if $A(u',z)\ge A(u'',z)$.
Thus the assertion is reduced to the case $\dim P\Gr_{w+i}^WH_{\R}=1$
or $2$.
In the two-dimensional case, take a basis $(v,\bar{v})$ of
$P\Gr_{w+i}^WH_{\C}$ which is compatible with the Hodge decomposition.
Then $av+a'\bar{v}$ is real if and only if $\bar{a}=a'$.
Using (2.3.1), the image of $u=\sum_{0\le k\le i}N^{i-k}u_k\in G_{\R}$
is then expressed as
$$\eqalign{&(d^kf(z)/dz^k,\,d^k\f(z)/dz^k)_{0\le k\le i},\cr
\hbox{if}\quad&f(z)=\msum_{0\le k\le i}\,a_kz^k/k!\quad\hbox{with}
\quad a_kv+\bar{a}_k\bar{v}=u_k,\cr}$$
where $\f(z)=\overline{f(\z)}$.
So the assertion is reduced to the following (where we assume
$f=\f$ and $n_1=n_2$ in the case $\dim P\Gr_{w+j}^WH_{\R}=1$).

\medskip\noindent
{\bf 2.5~Lemma.} {\it
Let $n, n_1, n_2$ be nonnegative integers such that $n_1+n_2>n$.
Let $\C [z]^{\le n}$ be the set of polynomials of degree $\le n$.
For $f\in\C[z]$ and $z_0\in\C$, set $y_0:=\Im\,z_0$, and
$$A(f,z_0):=\msum_{k\ge 0}\,|(d^kf/dz^k)(0)|y_0^k.$$
Then for any $a,a'\in\C$, there exist $C_{a,a'},\e>0$, $r>1$
such that for any $f(z)\in\C[z]^{\le n}$ and $z_0\in I_{r}$, we have
$$A(f,z_0)\le C_{a,a'},$$ if the following conditions are satisfied.}
$$\eqalign{|(d^k(f-a)/dz^k)(z_0)|&\le\e A(f,z_0)/y_0^k\quad
\hbox{for any}\,\,\,0\le k<n_1,\cr
|(d^k(\f-a')/dz^k)(z_0)|&\le\e A(f,z_0)/y_0^k\quad
\hbox{for any}\,\,\,0\le k<n_2.}
\leqno(2.5.1)$$

\medskip\noindent{\it Proof.}
Note first that $y_0\le|z_0|\le 2y_0$ since $r>1$.
We may assume $n_1\ge n_2$ exchanging $f$ and $\f$ if necessary.
For $f\in\C[z]$, set
$$b_k=(d^k(f-a)/dz^k)(z_0).$$
There is $g(z)\in\C [z]^{\le n-n_{1}}$ such that
$$f(z)-a=\msum_{0\le k<n_{1}\,}(b_k/k!)(z-z_0)^k+g(z)(z-z_0)^{n_{1}}.
$$
Using the binomial coefficients, we see that
$$A((z-z_0)^k,z_0)=A((z-\z_0)^k,z_0)\le C(k)y_0^k,\leqno(2.5.2)$$
where $C(k)>0$ depends only on $k$.
Applying the Leibniz rule to $g(z)(z-z_0)^{n_{1}}$ and using (2.5.1),
we then get
$$A(f-a,z_0)\le\e C_1A(f,z_0)+C_2A(g,z_0)y_0^{n_1},$$
where $C_1,C_2>0$ (depending only on $n,n_1$).
This implies
$$(1-\e C_1)A(f,z_0)\le C_2A(g,z_0)y_0^{n_1}+|a|.
\leqno(2.5.3)$$

On the other hand, the above equation can be modified as
$$\g(z)=(z-\z_0)^{-n_1}((\f(z)-a')+(a'-\bar{a}))-\msum_{0\le k<n_1\,}
(\bar{b}_k/k!)( z-\z_0)^{k-n_1}.$$
Then (2.5.1) implies similarly
$$|(d^k\g/dz^k)(z_0)|\le\e C_{3}A(f,z_0)/y_0^{k+n_1}+
C_4(|a|+|a'|)/y_0^{k+n_1}\quad\hbox{for}\,\,\,0\le k<n_2,$$
where $C_3,C_4>0$ (depending only on $n,n_1$).
Assuming $\e\le(2C_1)^{-1}$ and using (2.5.3), we then get
$$|(d^k\g/dz^k)(z_0)|\le\e C_5A(g,z_0)/y_0^k+C_6(|a|+|a'|)/
y_0^{k+n_1}\quad\hbox{for}\,\,\,0\le k\le n-n_1,$$
where $C_5,C_6>0$ (depending only on $n,n_1$).
Combined with (2.5.2), this implies
$$A(g,z_0)\le\e C_7A(g,z_0)+C_8(|a|+|a'|)/y_0^{n_1}
\quad\hbox{for}\,\,\,0\le k\le n-n_1,$$
where $C_7,C_8>0$ (depending only on $n,n_1$).
Assuming $\e\le(2C_7)^{-1}$, we deduce from this and (2.5.3)
$$A(f,z_0)\le C_9(|a|+|a'|),$$
where $C_9>0$ (depending only on $n,n_1$).
So the assertion follows.
This completes the proofs of Lemma~(2.5), Proposition~(2.4) and
Theorem~(2.1)(i).

\medskip\noindent
{\bf 2.6.~Proof of Theorem~(2.1)(ii).}
We have $v_k=0$ since $k>0$ and $v\in\Ker\,N$.
So it is enough to show
$$|\phi(\i(u);z,t')|_{k}\le\e A(u,z)/y^{k}\quad 
\hbox{if}\,\,\,{u}_{k'}^{(j)} = 0\,\,\,\hbox{for}\,\,\,k'\ge k.$$
Then the assertion is reduced to the nilpotent orbit case by (2.3.2)
(replacing $\e,\e',r$ if necessary).
In this case we have
$$|\phi(\i(u);z,t')|_{k} = 0\quad\hbox{if}\,\,\,u_{k'}^{(j)}=0\,\,\,
\hbox{for}\,\,\,k'\ge k,$$
because $\i_{i,j,k}=0$ for $k<0$ by (1.6.2).
So the assertion follows.

\medskip
From Theorem~(2.1) we can deduce

\medskip\noindent
{\bf 2.7.~Corollary.} {\it For any $p\in\VV_D^{\Inv}$, there is a
sufficiently small open neighborhood $U_p$ of $p$ in $\VV$ such that
$$U_p\cap\GG\subset\GG(p),\leqno(2.7.1)$$
where $\GG(p)$ is the section of $\GG$ passing through $p$ if
$p\in\GG$ and $\GG(p)=\emptyset$ otherwise.}

\medskip\noindent
{\it Proof.}
In case $p\in\GG$ we may assume $p=0$ using the action of $\GG(p)$ on
$\VV$.
Take $v\in G_{\C}^{<0}\cap\Ker\,N$ corresponding to $p$, and
apply Theorem~(2.1). Set
$$U_p=\big\{(v',t)\in G_{\C}^{<0}\times\D^m\,\big|\,|v'-v|<\e,
|t_1|<e^{2\pi r},|t_j|<\h\,(j>1)\big\},$$
where $\e,\h,r$ are as in Theorem~(2.1).
Then, for any $u\in G_{\Z}\setminus\Ker\,K$, (2.1.1) is satisfied
for some $k\ge 0$.
So the assertion follows (since $r>1$).

\medskip\noindent
{\bf 2.8.~Proof of Theorem~2.}
For the first assertion we have to show that for any $p\in\VV_D$ in
the notation of (1.2), there is an open neighborhood $U_p$ in $\VV$
such that
$$(U_p-U_p)\cap\GG\subset\GG(0),$$
where $\GG(0)$ denotes the zero section.
So it is enough to show that for $0\in\GG(0)$ over $0\in\S$,
there is an open neighborhood $U_0$ of $0$ in $\VV$ such that
$$U_0\cap(\GG\setminus\GG(0))=\emptyset.\leqno(2.8.1)$$
Then the assertion follows from Corollary~(2.7).

For the second assertion we have to show that for any
$p_1,p_2\in\VV_D^{\Inv}$ such that $\pi(p_1)=\pi(p_2)$ and
$p:=p_1-p_2\notin\GG$, there are open neighborhoods $U_i$ of $p_i$
in $\VV$ such that
$$(U_1-U_2)\cap\GG=\emptyset\quad\hbox{if}\,\,\,p\notin\GG,$$
where $\GG(p)\subset\GG$ is the section over $\S$ passing $p$
(shrinking $\S$ if necessary).
So it is enough to show that for any $p\in\VV_D^{\Inv}\setminus\GG$,
there is an open neighborhood $U_p$ of $p$ in $\VV$ such that
$$U_p\cap\GG=\emptyset.\leqno(2.8.2)$$
Then the assertion follows from Corollary~(2.7).

\medskip\noindent
{\bf 2.9.~Proof of Theorem~1.}
Taking a finite covering, the assertion is reduced to the case
where the local monodromy is unipotent.
Let $\n$ be an admissible normal function such that the local
cohomological invariant $\gamma(\n)$ is torsion, see (1.3).
Then there is a normal function $\m$ whose associated exact sequence
(1.3.1) splits after tensoring with $\Q$ and such that
$$\gamma(\n)=\gamma(\m).$$
Let $\nn,\mm$ be multivalued holomorphic sections of $\V$ lifting
$\n$ and $\m$ respectively.
Since $\m$ vanishes by the scaler extension $\Z\to\Q$, there is a
positive integer $m$ such that $m\1\m=0$, i.e. $\mm$ is a
multivalued section of ${1\over m}\G$. Set
$$\n'=\n-\m,\quad\nn'=\nn-\mm.$$
For the study of the zero locus of $\n$, it is then sufficient to
consider the intersection
$$\nn'\cap(\G-\mm)\,\,\,(\subset\nn'\cap\hbox{$1\over m$}\G),
\leqno(2.9.1)$$ instead of $\nn\cap\G$.
Here a section is identified with its image.
Since $\gamma(\n')=0$, we may assume that $\nn'$ is univalent and is
extended to a section of $\GG$.
As explained in (1.3), we have by a recent theorem of Green,
Griffiths and Kerr [8]
$$\nn'(0)\in J^{\GGK}_{\S}(\H)^0.$$
Then, applying Corollary~(2.7) to $\nn'(0)$ and
${1\over m}\GG\subset\VV$, there is a sufficiently small open
neighborhood $U_{\nn'(0)}$ of $\nn'(0)$ such that
$$\hbox{${1\over m}$}\GG\cap U_{\nn'(0)}=
\cases{\mm'\cap U_{\nn'(0)}& if $\,\,\nn'(0)\in{1\over m}\GG$,\cr
\emptyset & if $\,\,\nn'(0)\notin{1\over m}\GG$,}
\leqno(2.9.2)$$
where $\mm'$ is the section of ${1\over m}\GG$ passing through
$\nn'(0)\in{1\over m}\GG$.

If $\nn'(0)\notin{1\over m}\GG$ or $\nn'(0)\in{1\over m}\GG$ with
$\gamma(\n)\ne 0$, then (2.9.1) is empty (shrinking $\S$ if necessary)
since $\mm'$ is not a section of $\G-\mm$ in the last case.
(Indeed, if $\mm'+\mm$ is a multivalued section of $\G$, then $\mm'$
coincides with $-\mm$ replacing $\mm$ if necessary, and hence $\mm$
is univalent. But this contradicts the condition that
$\gamma(\m)=\gamma(\n)\ne 0$.)
This implies that the zero locus of $\n$ is empty on a neighborhood
of $0$.

If $\nn'(0)\in{1\over m}\GG$ with $\gamma(\n)=0$, then we may assume
$\mm=0$, $\nn'=\nn$, $m=1$, and furthermore $\nn(0)=0$ replacing
$\nn$ if necessary. So the assertion follows from (2.9.2).

This completes the proof of Theorem~1.

\medskip\noindent
{\bf 2.11.~Proof of Theorem~3.}
By definition [8] (see also [10]), the N\'eron model
$J_{\S}^{\GGK}(\H)$ is the union of
$Y_{\nu}:=J_{\S}^{\GGK}(\H)^0|_{U(\nu)}$,
where $\nu$ runs over normal functions defined over open subsets
$U(\nu)\setminus D$ which are admissible with respect to $U(\nu)$.
We have the identifications
$$Y_{\nu}|_{U(\nu)\cap U(\nu')\setminus D}\cong
Y_{\nu'}|_{U(\nu)\cap U(\nu')\setminus D},$$
induced by $\nu-\nu'$ over $U(\nu)\cap U(\nu')\setminus D$ if the
cohomology classes of $\nu$ and $\nu'$ are different.
In case they coincide, then the isomorphisms is extended over
$U(\nu)\cap U(\nu')$.
Moreover we may assume that the $\nu$ are torsion since $w=-1$.
So it is enough to show that for any $p_1,p_2\in\VV_D^{\Inv}$
with $\pi(p_1)=\pi(p_2)$, there exists respectively open
neighborhoods $U_1,U_2$ of $p_1,p_2$ in $\VV$ such that
$$(U_1-U_2)\cap\hbox{${1\over m}$}\GG$$
is contained in the sections of ${1\over m}\GG$ {\it over} $\S$
shrinking $\S$ if necessary.
Then the assertion follows from Corollary~(2.7).

\medskip\noindent
{\bf 2.12.~Remark.}
The norm in (2.1.1) uses the trivialization (1.8.2),
and is quite different from the one in [11] and [4] using the Weil
operator at each point of $S$.
It does not seem easy to generalize the method in this paper to
the case $D$ is a divisor with normal crossings since
at least the multi $SL_2$-orbit theorem [4] will be needed.
Note that this case is treated in [2].

\bigskip\bigskip
\centerline{{\bf References}}

\medskip
{\mfont
\item{[1]}
P.~Brosnan and G.~J.~Pearlstein,
The zero locus of an admissible normal function (preprint
math.AG/0604345), to appear in Ann. of Math.

\item{[2]}
P.~Brosnan and G.~J.~Pearlstein,
Zero loci of admissible normal functions with torsion singularities
(preprint arXiv:0803.3365).

\item{[3]}
J.~Carlson, Extensions of mixed Hodge structures, in Journ\'ees de
G\'eom\'etrie Alg\'ebrique d'Angers 1979, Sijthoff-Noordhoff Alphen
a/d Rijn, 1980, pp. 107--128.

\item{[4]}
E.~Cattani, A.~Kaplan, W.~Schmid, Degeneration of Hodge structures,
Ann. of Math., 123 (1986), 457--535.

\item{[5]}
P.~Deligne, Equations diff\'erentielles\`a points singuliers
r\'eguliers, Lect. Notes in Math. vol. 163, Springer, Berlin, 1970.

\item{[6]}
P.~Deligne, Th\'eorie de Hodge II, Publ. Math. IHES 40 (1971),
5--58.

\item{[7]}
F.~El Zein and S.~Zucker, Extendability of normal functions associated
to algebraic cycles, in Topics in transcendental algebraic geometry,
Ann. Math. Stud., 106, Princeton Univ. Press, Princeton, N.J., 1984,
pp. 269--288.

\item{[8]}
M.~Green, P.~Griffiths and M.~Kerr,
N\'eron models and limits of Abel-Jacobi mappings (preprint).

\item{[9]}
M.~Kashiwara, A study of variation of mixed Hodge structure,
Publ.\ RIMS,\ Kyoto Univ. 22 (1986), 991--1024.

\item{[10]}
M.~Saito,
Admissible normal functions, J.\ Algebraic Geom. 5 (1996),
235--276.

\item{[11]}
W.~Schmid, Variation of Hodge structure: The singularities of the
period mapping, Inv. Math. 22 (1973), 211--319.

\item{[12]}
J.H.M.~Steenbrink, Limits of Hodge structures, Inv. Math. 31 (1975/76), 
229--257.

\item{[13]}
J.H.M.~Steenbrink and S.~Zucker, Variation of mixed Hodge structure I,
Inv.\ Math. 80 (1985), 489--542.

\item{[14]}
S.~Zucker, Generalized intermediate Jacobians and the theorem on
normal functions, Inv. Math. 33 (1976),185--222.

\medskip
RIMS Kyoto University, Kyoto 606-8502 Japan

\medskip
\vers
}
\bye